\documentclass{amsproc}%
\usepackage{amssymb}
\usepackage{graphicx}
\usepackage{amscd}
\usepackage{amsmath}
\usepackage{amsfonts}%
\theoremstyle{plain}

\numberwithin{equation}{section}

\begin{document}
\title[Which self-maps appear as lattice anti-endomorphisms?]{Which self-maps appear as lattice anti-endomorphisms?}
\author{Stephan Foldes}
\address{Institute of Mathematics, Tampere University of Technology, PL 553, Tampere
33101, Finland}
\email{stephan.foldes@tut.fi}
\author{Jen\H{o} Szigeti}
\address{Institute of Mathematics, University of Miskolc, Miskolc, Hungary 3515}
\email{matjeno@uni-miskolc.hu}
\thanks{The second author was supported by OTKA K101515 and his research was carried
out as part of the TAMOP-4.2.1.B-10/2/KONV-2010-0001 project with support by
the European Union, co-financed by the European Social Fund.}

\begin{abstract}
Let $f:A\rightarrow A$ be a self-map of the set $A$. We give a necessary and
sufficient condition for the existence of a lattice structure $(A,\vee
,\wedge)$ on $A$ such that $f$ becomes a lattice anti-endomorphism with
respect to this structure.

\end{abstract}
\maketitle

\noindent1. INTRODUCTION

\bigskip

A \textit{partially ordered set (poset)} is a set $P$\ together with a
reflexive, antisymmetric, and transitive (binary) relation $r\subseteq P\times
P$. For $(x,y)\in r$ we write $x\leq_{r}y$ or simply $x\leq y$. If $r\subseteq
r^{\prime}$ for the partial orders $r$ and $r^{\prime}$ on $P$, then
$r^{\prime}$ is an \textit{extension} of $r$. A map $p:P\longrightarrow P$ is
\textit{order-preserving (order reversing)} if $x\leq y$ implies $p(x)\leq
p(y)$ ($p(y)\leq p(x)$ respectively) for all $x,y\in P$. The poset $(P,\leq)$
is a \textit{lattice} if any two elements $x,y\in P$ have a unique least upper
bound (lub) $x\vee y$ and a unique greatest lower bound (glb) $x\wedge y$ (in
$P$).

In the present paper we consider a self-map $f:A\longrightarrow A$\ of a set
$A$. A list $x_{1},\ldots,x_{n}$ of distinct elements from $A$ is a
\textit{cycle} (of \textit{length} $n$) with respect to $f$ if $f(x_{i}%
)=x_{i+1}$ for each $1\leq i\leq n-1$ and also $f(x_{n})=x_{1}$. A
\textit{fixed point} of the function $f$ is a cycle of length $1$, i.e. an
element $x_{1}\in A$ with $f(x_{1})=x_{1}$. A cycle that is not a fixed point
is \textit{proper}.

If $(A,\vee,\wedge)$ is a lattice (on the set $A$) such that $f(x\vee
y)=f(x)\wedge f(y)$ and $f(x\wedge y)=f(x)\vee f(y)$ for all $x,y\in A$, then
$f$ is a \textit{lattice anti-endomorphism} of $(A,\vee,\wedge)$. The square
$f^{2}=f\circ f$ of a lattice anti-endomorphism is an ordinary lattice
endomorphism of $(A,\vee,\wedge)$. A lattice anti-endomorphism is an
order-reversing map (with respect to the order relation of the lattice), but
the converse is not true in general. More details about lattices and lattice
(anti)endomorphisms can be found in [G].

For a proper cycle $x_{1},\ldots,x_{n}\in A$ with respect to a lattice
anti-endomorphism $f$, if we put%
\[
p=x_{1}\vee x_{2}\vee\cdots\vee x_{n}\text{ and }q=x_{1}\wedge x_{2}%
\wedge\cdots\wedge x_{n},
\]
then $p\neq q$. The equalities%
\[
f(p)=f(x_{1})\wedge f(x_{2})\wedge\cdots\wedge f(x_{n})=x_{2}\wedge
\cdots\wedge x_{n}\wedge x_{1}=q,
\]%
\[
f(q)=f(x_{1})\vee f(x_{2})\vee\cdots\vee f(x_{n})=x_{2}\vee\cdots\vee
x_{n}\vee x_{1}=p
\]
show that $p$ and $q$ form a cycle of length $2$\ with respect to $f$. If
$u,v\in A$ are distinct fixed points of $f$, then $p=u\vee v$ and $q=u\wedge
v$ form again a cycle of length $2$\ with respect to $f$:%
\[
f(p)=f(u)\wedge f(v)=u\wedge v=q\text{ and }f(q)=f(u)\vee f(v)=u\vee v=p.
\]
It follows that any lattice anti-endomorphism having a proper cycle or having
at least two fixed points must have a cycle of length $2$.

We prove that the above combinatorial property completely characterizes the
possible lattice anti-endomorphisms. More precisely, for a map
$f:A\longrightarrow A$\ there exists a lattice $(A,\vee,\wedge)$ on $A$ such
that $f$ is a lattice anti-endomorphism of $(A,\vee,\wedge)$ if and only if
$f$ has a cycle of length $2$ or $f$ has no proper cycles and has at most one
fixed point.

In section 2 we give a relatively short and self-contained proof for the
mentioned characterization of lattice anti-endomorphisms by using conditional
lattices. First we prove that every self-map is a conditional lattice
anti-endomorphism with respect to some conditional lattice structure on the
base set.

Section 3 is devoted to the study of lattice-order extensions preserving the
anti-monotonicity of a given self-map $f:A\longrightarrow A$. We show that
certain partial orders on $A$\ can be extended to lattice-orders making $f$ an
anti-endomorphism. Thus we obtain a broad generalization of the pure existence
result in section 2. Our treatment in section 3 follows the lines of [Sz2],
where a combinatorial characterization of ordinary lattice endomorphisms can
be found. The results of [FSz] served as a starting point in [Sz2], while in
section 3 we build on [Sz1]. The construction in the proof of Theorem 3.8 is
based on the use of the maximal antimonotonicity preserving (AMP) extensions
of $r$\ in a partially anti-ordered unary algebra $(A,f,\leq_{r})$. Such
extensions were completely determined in [L] and [Sz1]. In order to make the
exposition more self-contained, we present the necessary background about the
mentioned AMP extensions.

\bigskip

\noindent2. CONDITIONAL\ LATTICES AND\ ANTI-ENDOMORPHISMS

\bigskip

Recall that a partial order $(A,\leq_{r})$\ is called a \textit{conditional
lattice} if for any two elements $x$ and $y$ the following holds: if the two
elements have a common upper bound, then they have a least upper bound,
denoted $x\vee y$, and if they have a common lower bound, then they have a
greatest lower bound, denoted $x\wedge y$. In those cases we also say that
$x\vee y$ or $x\wedge y$ \textit{exist. }All lattices are conditional
lattices, so are all antichains, also any poset that is a union of disjoint
chains with no comparabilities between elements from the different chains
(unordered sum of chains).

An \textit{anti-endomorphism} of a conditional lattice $(A,\leq_{r})$ is a
self-map $f:A\longrightarrow A$ such that for all $x,y\in A$ the following
hold:%
\[
x\vee y\text{ exists }\Longleftrightarrow\text{ }f(x)\wedge f(y)\text{ exists
and equals }f(x\vee y)
\]%
\[
x\wedge y\text{ exists }\Longleftrightarrow\text{ }f(x)\vee f(y)\text{ exists
and equals }f(x\wedge y)
\]

In what follows we use the language of directed graphs. The \textit{digraph of
a self-map} $f:A\longrightarrow A$ has node set $A$ and arrow set
$\{(x,f(x))\mid x\in A\}$.

\bigskip

\noindent\textbf{2.1. Conditional Lattice Lemma.}\textit{ Every function
}$f:A\longrightarrow A$\textit{ is an anti-endomorphism of some conditional
lattice on }$A$\textit{.}

\bigskip

\noindent\textbf{Proof.\ }Without loss of generality assume connectedness,
i.e. that the digraph $D$ on node set $A$ with arrows $(x,f(x))$ is simply
connected. The \textit{tree components} are the connected components obtained
from $D$ by removing the arrows of any (unique if existing) directed cycle.
For any tree component $T$ the \textit{sons} of a node $y$ are the nodes $x$
in $T$ with $f(x)=y$. For each node $x$ of $T$ choose a linear order $\leq
_{x}$ on the set of its sons, having a first son and a last son, if $x$ has
any son. Let $P_{x}$ denote the maximal directed path in $T$ starting from
$x$. For nodes $x$, $y$ in $T$ let $z$ be the first common node of $P_{x}$ and
$P_{y}$, and write $x\leq_{\mathrm{lex}}y$ whenever (for distance denoted by
$d$)%
\[
x=y\text{ or }d(x,z)<d(y,z)
\]
or%
\[
x\neq y\text{, }d(x,z)=d(y,z)\text{ and }x^{\prime}\leq_{z}y^{\prime}%
\]%
\[
\text{for the nodes }x^{\prime}\text{, }y^{\prime}\text{ preceeding }z\text{
immediately in }P_{x}\text{, }P_{y}\text{.}%
\]
This defines a linear order on the nodes of $T$, call \textit{lexicographic}
(lex) order. Define partial orders on $A$, distinguishing two cases.

\textit{Case 1}: $D$ has no cycle. Let $g:A\longrightarrow\mathbb{Z}$ be a
grading of $A$, i.e. $g(f(x))=g(x)+1$ for all $x\in A$. Let $N(i)=g^{-1}(i)$
for $i\in\mathbb{Z}$. Let $E$ and $O$ be the set of nodes of even and odd
grade, respectively. Order $E$\ lexicographically $O$
reverse-lexicographically, and let all elements of $E$ be less than all
elements of $O$. This is a linear order on $A$, $f$ is order-reversing.

\textit{Case 2}: $D$ has a cycle with $m$ distinct nodes $a_{1},\ldots,a_{m}$
indexed by mod $m$ integers, $m\geq1$, such that $f(a_{i})=a_{i+1}$. There is
a tree component $T_{i}$ for each $a_{i}$. Let $N(i,k)$ be the set of nodes of
$T_{i}$ at distance $k\geq0$ form $a_{i}$. Order $N(i,k)$ lex or reverse lex
according to wether $k$ is even or odd.

Let all elements of $N(i,k)$ be less (greater) than all elements of
$N(i+1,k+1)$ if $k$ is even (odd), with no further comparabilities.

This order on $A$ is a conditional lattice, $f$ is an anti-endomorphism.
$\square$

\bigskip

\noindent\textbf{2.2. Theorem.}\textit{ A function }$f:A\longrightarrow
A$\textit{\ is a lattice anti-endomorphism of some lattice structure on }%
$A$\textit{ if and only if either }$f$\textit{\ transposes some pair of
elements or it does not induce a permutation on any finite set of at least two
elements.}

\bigskip

\noindent\textbf{Proof.}

\textit{Case 1}: $f$ has a single fixed point $a$ in a connected component $K$
of its digraph $D$, and the other connected components $K_{i}$, $i\in\sigma$,
indexed by some ordinal $\sigma$, have no cycle. Removing the loop arrow from
$K$ we obtain a tree whose node set can be ordered lexicographically as in the
proof of the Lemma. Classify the nodes of $K$ into even and odd sets $E$, $O$
according to the parity of their distance from $a$. Classify the nodes of each
$K_{i}$ into even and odd sets $E_{i}$, $O_{i}$ according to a $\mathbb{Z}%
$-grading as in the proof of the Lemma.

Let $E$ and the $E_{i}$'s be ordered lexicographically, $O$ and the $O_{i}$'s
reverse lexicographically. Order $A$ by%
\[
\cdots<O_{2}<O_{1}<O_{0}<O<E<E_{0}<E_{1}<E_{2}<\cdots,
\]
i.e. the odd and even sets now appear as succeeding intervals of a linear
order on$~A$.

\textit{Case 2}: If $f$ has no cycle, not even a fixed point, then simply omit
$K$, $O$, $E$ in the above construction.

\textit{Case 3}: $f$ has a $2$-cycle, i.e. for some $a_{1}\neq a_{2}$,
$f(a_{1})=a_{2}$, $f(a_{2})=a_{1}$. Let $K$ be the connected component of the
digraph of $f$ containing this $2$-cycle, and let $B$ denote the elements of
$A$ not in $K$. Then both $B$ and $A\smallsetminus B$ are closed under $f$,
and by the above Lemma there is a conditional lattice structure $L$ on $B$, on
which $f$ induces an anti-endomorphism. Let $T_{1}$, $T_{2}$ be the trees
obtained from $K$ by removing the arrows of the $2$-cycle between $a_{1}$ and
$a_{2}$. For $i=1,2$ and $k\geq1$ let $N(i,k)$ denote the set of nodes of
$T_{i}$ at distance $k$ form $a_{i}$. Fix a lexicographic order on the nodes
of each $T_{i}$ as in the proof of the Lemma. Order $N(i,k)$ lexicographically
for $k$ even, reverse lexicographically for $k$ odd. Order $A$ by%
\[
\underset{N(1,k)\text{ }k\text{ odd decreasing}}{\underbrace{\cdots
<N(1,5)<N(1,3)<N(1,1)}}<\underset{N(2,k)\text{ }k\text{ odd decreasing}%
}{\underbrace{\cdots<N(2,5)<N(2,3)<N(2,1)}}<a_{1}<L
\]%
\[
L<a_{2}<\underset{N(2,k)\text{ }k\text{ even increasing}}{\underbrace
{N(2,2)<N(2,4)<N(2,6)<\cdots}}\underset{N(1,k)\text{ }k\text{ even
increasing}}{<\underbrace{N(1,2)<N(1,4)<N(1,6)<\cdots}}%
\]
i.e. a chain made up by succeeding intervals $N(i,k)$, $k$ odd, followed by
the element $a_{1}$, the conditional lattice $L$, then the element $a_{2}$,
then by a chain made up by succeeding intervals $N(i,k)$, $k$ even. This order
is a lattice, of which $f$ is an anti-endomorphism. $\square$

\bigskip

\noindent3. LATTICE\ ORDER\ EXTENSIONS\ AND ANTI$-$ENDOMORPHISMS

\bigskip

The following definitions appear in [FSz] and [Sz1]. The treatment in the
paper [JPR] also starts with similar considerations.

Let $f:A\longrightarrow A$ be a function and define the equivalence relation
$\sim_{f}$ as follows: for $x,y\in A$ let $x\sim_{f}y$ if $f^{k}(x)=f^{l}(y)$
for some integers $k\geq0$ and $l\geq0$ with $f^{0}$ meaning the identity
function. The equivalence class $\left[  x\right]  _{f}$\ of an element $x\in
A$ is called the $f$\textit{-component} of $x$. We note that $f(\left[
x\right]  _{f})\subseteq\left[  x\right]  _{f}$.

An element $c\in A$ is called \textit{cyclic} with respect to $f$ (or
$f$-cyclic), if $f^{m}(c)=c$ for some integer $m\geq1$. The \textit{period} of
a cyclic element $c$ is%
\[
n_{f}(c)=\min\{m\mid m\geq1\text{ and }f^{m}(c)=c\}.
\]
Obviously, $f^{k}(c)=f^{l}(c)$ holds if and only if $k-l$ is divisible by
$n_{f}(c)$. A cyclic element $c$ or the $f$\textit{-cycle} $\{c,f(c),\ldots
,f^{n_{f}(c)-1}(c)\}$ is called \textit{proper} if $n_{f}(c)\geq2$. The
$f$-orbit $\{x,f(x),\ldots,f^{k}(x),\ldots\}$ of $x$ is finite if and only if
$\left[  x\right]  _{f}$\ contains a cyclic element. If $c_{1},c_{2}\in\left[
x\right]  _{f}$ are cyclic elements, then $n_{f}(c_{1})=n_{f}(c_{2})=n_{f}(x)$
and this number is called the \textit{period }of $x$. If the $f$-orbit of $x$
is infinite, then define $n_{f}(x)=\infty$. Clearly, $x\sim_{f}y$ implies that
$n_{f}(x)=n_{f}(y)$. If $f$ has a cycle of odd period, then it also defines a
cycle of $f^{2}=f\circ f$ with the same elements and the same (odd) period. If
$f$ has a cycle of even period $n=2m\geq2$, then it is a disjoint union of two
$f^{2}$-cycles of period $m$. If $n_{f}(x)$ is odd, then $\left[  x\right]
_{f}=\left[  x\right]  _{f^{2}}=\left[  f(x)\right]  _{f^{2}}$ and $n_{f^{2}%
}(x)=n_{f}(x)$. If $n_{f}(x)$ is even or $n_{f}(x)=\infty$, then
$x\nsim_{f^{2}}f(x)$, $\left[  x\right]  _{f}=\left[  x\right]  _{f^{2}}%
\cup\left[  f(x)\right]  _{f^{2}}$ and $n_{f^{2}}(x)=\frac{1}{2}n_{f}(x)$.
Notice that $f(\left[  x\right]  _{f^{2}})\subseteq\left[  f(x)\right]
_{f^{2}}$.

The set of $f$\textit{-incomparable pairs} is $\pi=\alpha\cup\beta$, where%
\[
\alpha=\{(x,y)\in A\times A\mid x\nsim_{f}y\text{, }n_{f}(x)\text{ and }%
n_{f}(y)\text{ are odd integers}\},
\]%
\[
\beta=\{(x,y)\in A\times A\mid f^{2k+m}(x)=f^{m}(y)\neq f^{m}(x)=f^{2l+m}%
(y)\text{ for some }k,l,m\geq0\}.
\]
If $(x,y)$ is not $f$-incomparable (i.e. $(x,y)\notin\pi$), then we say that
$(x,y)$ is $f$\textit{-comparable}. The following properties can easily be checked.

\noindent1. If $(x,y)$\ is $f$-incomparable, then $(y,x)$\ is also $f$-incomparable.

\noindent2. If $(f(x),f(y))$\ is $f$-incomparable, then $(x,y)$\ is also $f$-incomparable.

\noindent3. If $f$ has more than one fixed point, then $\alpha\neq\varnothing$.

\noindent4. If $(x,y)\in\beta$ and $m$ is as in the definition of $\beta$,
then $f^{m}(x)$ and $f^{m}(y)$ are

\noindent\ \ \ \ (different) elements of the same $f^{2}$-cycle.

\noindent5. If $f^{2}$ has a proper cycle, then $\beta\neq\varnothing$.

\noindent6. If $f$ has at most one fixed point and $f^{2}$ has no proper
cycle, then $\alpha=\beta=\varnothing$

\noindent\ \ \ \ (hence $\pi=\varnothing$).

\bigskip

\noindent\textbf{3.1. Lemma.}\textit{ For }$x,y\in A$\textit{ the following
conditions are equivalent}

1.\textit{ }$(x,y)\in\beta$\textit{,}

2.\textit{ }$(f(x),f(y))\in\beta$\textit{,}

3.\textit{ }$3\leq n_{f}(x)\neq\infty$\textit{, }$x\sim_{f^{2}}y$\textit{ and
}$f^{t}(x)\neq f^{t}(y)$\textit{ for all integers }$t\geq0$\textit{.}

\bigskip

\noindent\textbf{Proof.} (1)$\Longrightarrow$(2)\&(3). Now $f^{2k+m}%
(x)=f^{m}(y)\neq f^{m}(x)=f^{2l+m}(y)$ for some $k,l,m\geq0$. It follows that
$k,l\geq1$. Since $f^{2k+m}(x)=f^{m}(y)$ implies $f^{2k+m+1}(x)=f^{m+1}(y)$
and one of $m$ and $m+1$ is even, we obtain that $x\sim_{f^{2}}y$. Clearly,
$f^{m}(y)\neq f^{m}(x)$ implies that $f^{t}(x)\neq f^{t}(y)$ for all $0\leq
t\leq m$. In view of%
\[
f^{2k+2l-1}(f^{m+1}(x))\!=\!f^{2k+2l+m}(x)\!=\!f^{2l}(f^{2k+m}(x))\!=\!f^{2l}%
(f^{m}(y))\!=\!f^{2l+m}(y)\!=\!f^{m}(x)
\]
and%
\[
f^{2k+2l-1}(f^{m+1}(y))\!=\!f^{2k+2l+m}(y)\!=\!f^{2k}(f^{2l+m}(y))\!=\!f^{2k}%
(f^{m}(x))\!=\!f^{2k+m}(x)\!=\!f^{m}(y)
\]
the equality $f^{m+1}(x)=f^{m+1}(y)$ would imply $f^{m}(x)=f^{m}(y)$, a
contradiction. It follows that $f^{2k+m}(f(x))\!=\!f^{m}(f(y))\!\neq
\!f^{m}(f(x))\!=\!f^{2l+m}(f(y))$ and $(f(x),f(y))\!\in~\beta$. Starting from
$(f(x),f(y))\!\in\!\beta$ and replacing the triple $(x,y,m)$ by
$(f(x),f(y),m\!+\!1)$ in the above argument give that $f^{m+2}(x)\neq
f^{m+2}(y)$. Thus $f^{t}(x)\neq f^{t}(y)$ for all $t\geq m$. As a consequence
of $f^{2k+2l+m}(x)=f^{m}(x)$, we deduce that $f^{m}(x)$ is $f$-cyclic. Since
$f^{m+2}(x)=f^{m}(x)$ would imply $f^{m}(x)=f^{2k+m}(x)=f^{m}(y)$, we obtain
that $3\leq n_{f}(x)\neq\infty$.

(2)$\Longrightarrow$(1) is straightforward.

(3)$\Longrightarrow$(1). $x\sim_{f^{2}}y$ imply the existence of integers
$r,s\geq0$ such that $f^{2r}(x)=f^{2s}(y)$. Since $3\leq n_{f}(x)\neq\infty$,
there exists an integer $t\geq0$ such that $f^{2t+2r}(x)=f^{2t+2s}(y)$ is
$f$-cyclic of period $n_{f}(x)$. In view of $2t+2r\leq2t+2r+2s$ and
$2t+2s\leq2t+2r+2s$, the relations%
\[
f^{2t+2r+2s}(x)\sim_{f^{2}}x\sim_{f^{2}}y\sim_{f^{2}}f^{2t+2r+2s}(y)
\]
imply that $f^{2t+2r+2s}(x)$ and $f^{2t+2r+2s}(y)$ are different $f$-cyclic
elements in $\left[  x\right]  _{f^{2}}$. It follows that $f^{2t+2r+2s}(x)$
and $f^{2t+2r+2s}(y)$ are in the same $f^{2}$-cycle and%
\[
f^{2k}(f^{2t+2r+2s}(x))=f^{2t+2r+2s}(y)\neq f^{2t+2r+2s}(x)=f^{2l}%
(f^{2t+2r+2s}(y))
\]
for some $k,l\geq1$ ($k+l$ is a multiple of $n_{f^{2}}(x)=n_{f^{2}}(y)$). Thus
$(x,y)\in\beta$. $\square$

\bigskip

\noindent\textbf{3.2. Corollary.}\textit{ For any }$a\in A$\textit{ the
relation }$(\left[  a\right]  _{f^{2}}\times\left[  a\right]  _{f^{2}%
})\smallsetminus\beta$\textit{ is an equivalence relation on the set }$\left[
a\right]  _{f^{2}}$\textit{. If }$3\leq n_{f}(a)\neq\infty$\textit{, then
there are exactly }$n_{f^{2}}(a)$\textit{ equivalence classes with respect to
}$(\left[  a\right]  _{f^{2}}\times\left[  a\right]  _{f^{2}})\smallsetminus
\beta$\textit{.}

\bigskip

\noindent\textbf{Proof.} If $1\leq n_{f}(a)\leq2$ or $n_{f}(a)=\infty$, then
$(\left[  a\right]  _{f^{2}}\times\left[  a\right]  _{f^{2}})\cap
\beta=\varnothing$ and our claim holds. If $3\leq n_{f}(a)\neq\infty$, then%
\[
(\left[  a\right]  _{f^{2}}\times\left[  a\right]  _{f^{2}})\smallsetminus
\beta=\{(x,y)\mid x,y\in\left[  a\right]  _{f^{2}}\text{ and }f^{t}%
(x)=f^{t}(y)\text{ for some }t\geq0\},
\]
which is obviously reflexive, symmetric and transitive. It is straightforward
to see that $\{a,f^{2}(a),\ldots,f^{2n_{f^{2}}(a)}(a)\}$\ is a complete
irredundant system of representatives with respect to the equivalence relation
$(\left[  a\right]  _{f^{2}}\times\left[  a\right]  _{f^{2}})\smallsetminus
\beta$. $\square$

\bigskip

A\textit{ partially anti-ordered unary algebra} is a triple $(A,f,\leq_{r})$,
where $r$ is a partial order on $A$ and $f:A\longrightarrow A$\ is an order
reversing map with respect $r$. Important facts about the close relationship
between the cycle and the order structure in such (and similar) triples can be
found in [FSz],[JPR],[L] and [Sz1].

\bigskip

\noindent\textbf{3.3. Proposition }(see [Sz1])\textbf{.}\textit{\ If
}$(A,f,\leq_{r})$\textit{\ is a partially anti-ordered unary algebra and
}$(x,y)\in\pi$\textit{\ is an }$f$\textit{-incomparable pair, then }%
$x$\textit{ and }$y$\textit{\ are incomparable with respect to }$r$\textit{.}

\bigskip

\noindent\textbf{3.4. Corollary }(see [Sz1])\textbf{.}\textit{\ If }%
$(A,f,\leq_{r})$\textit{\ is a partially anti-ordered unary algebra,
}$r\subseteq R$\textit{ is an }AMP\textit{ extension of }$r$\textit{ and
}$(x,y)\in\pi$\textit{\ is an }$f$\textit{-incomparable pair, then }%
$x$\textit{ and }$y$\textit{\ are incomparable with respect to }$R$\textit{.}

\bigskip

An AMP extension $R$\ of $r$ is called $f$\textit{-maximal}, if $x\leq_{R}y$
or $y\leq_{R}x$ for all $f$-comparable pairs $(x,y)\in A\times A$ (that is
$(A\times A)\setminus\pi\subseteq R\cup R^{-1}$). Corollary 3.4 implies that
any $f$-maximal AMP extension of $r$ is maximal with respect to containment.
We shall make use of the following notations:%
\[
\mathcal{M}(A,f,\leq_{r})=\{R\mid r\subseteq R\text{ and }R\text{\ is
}f\text{-maximal AMP extension\ of }r\},
\]%
\[
\mathcal{L}(A,f,\leq_{r})=\{R\mid r\subseteq R\text{ and }R\text{\ is an AMP
linear order}\}.
\]
Obviously, $\mathcal{L}(A,f,\leq_{r})\subseteq\mathcal{M}(A,f,\leq_{r})$. If
the function $f^{2}$ has a proper cycle or $f$ has more than one fixed point,
then $\pi\neq\varnothing$ and Corollary 3.4 gives that $\mathcal{L}%
(A,f,\leq_{r})=\varnothing$. If $f^{2}$ has no proper cycle and $f$ has at
most one fixed point, then $\pi=\alpha\cup\beta=\varnothing$ implies that
$\mathcal{L}(A,f,\leq_{r})=\mathcal{M}(A,f,\leq_{r})$.

\bigskip

\noindent\textbf{3.5. Theorem }(see [L])\textbf{.}\textit{\ If }$f^{2}%
$\textit{ has no proper cycle and }$f$\textit{ has at most one fixed point,
then }$\mathcal{L}(A,f,\leq_{r})\neq\varnothing$\textit{.}

\bigskip

\noindent\textbf{3.6. Theorem }(see [Sz1])\textbf{.}\textit{\ If }%
$(A,f,\leq_{r})$\textit{\ is an arbitrary partially anti-ordered unary
algebra, then }$\mathcal{M}(A,f,\leq_{r})\neq\varnothing$\textit{ and the
elements of }$\mathcal{M}(A,f,\leq_{r})$\textit{ are exactly the maximal (with
respect to containment) elements in the set of all }AMP\textit{ extensions of
}$r$\textit{.}

\bigskip

\noindent\textbf{3.7. Theorem.}\textit{\ Let }$f:A\longrightarrow
A$\textit{\ be a function such that }$f^{2}$\textit{ has no proper cycles and
}$f$\textit{\ has at most one fixed point. Then there exists a distributive
lattice }$(A,\vee,\wedge)$\textit{ on }$A$\textit{ such that }$f$\textit{ is a
lattice anti-endomorphism of }$(A,\vee,\wedge)$\textit{.}

\bigskip

\noindent\textbf{Proof.} A linearly ordered set is a distributive lattice and
an order reversing map with respect to this linear order is a lattice
anti-endomorphism. Thus the existence of the desired distributive lattice
(chain) is an immediate consequence of Lengv\'{a}rszky's Theorem 3.5.
$\square$

\bigskip

\noindent\textbf{3.8. Theorem.}\textit{\ Let }$(A,f,\leq_{r})$\textit{\ be a
partially anti-ordered unary algebra such that }$f$\textit{\ has a cycle
}$\{p,q\}$\textit{\ of length }$2$\textit{. If }$x$\textit{ and }$y$\textit{
are }$r$\textit{-incomparable for all }$x,y\in A$\textit{ with }$\left[
x\right]  _{f^{2}}\neq\left[  y\right]  _{f^{2}}$\textit{ and }$2\neq
n_{f}(x)\neq\infty$\textit{, then there exists an extension }$R$\textit{ of
}$r$\textit{ such that }$(A,\leq_{R})$\textit{ is a lattice and }$f$\textit{
is a lattice anti-endomorphism of }$(A,\leq_{R})$\textit{.}

\bigskip

\noindent\textbf{Proof.} Let%
\[
A_{0}=\{x\in A:2\neq n_{f}(x)\neq\infty\}
\]
and%
\[
A_{\ast}=A\setminus A_{0}=\{x\in A:n_{f}(x)=2\text{ or }n_{f}(x)=\infty\}.
\]
We have either $\left[  x\right]  _{f^{2}}\subseteq\left[  x\right]
_{f}\subseteq A_{0}$ or $\left[  x\right]  _{f^{2}}\subseteq\left[  x\right]
_{f}\subseteq A_{\ast}$ for all $x\in A$. Clearly, both $A_{0}$ and $A_{\ast}$
are closed with respect to the action of $f$, i.e. $f(A_{0})\subseteq A_{0}$
and $f(A_{\ast})\subseteq A_{\ast}$.

Take an arbitrary $f$-maximal AMP extension $R$ of $r$ (Theorem 3.6 ensures
the existence of such $R$). Since $\pi\cap(A_{\ast}\times A_{\ast
})=\varnothing$ implies%
\[
A_{\ast}\times A_{\ast}\subseteq(A\times A)\setminus\pi\subseteq R\cup
R^{-1},
\]
we deduce that $R_{\ast}=R\cap(A_{\ast}\times A_{\ast})$ is a linear order on
$A_{\ast}$. In view of $p,q\in A_{\ast}$, we may assume $p\leq_{R_{\ast}}q$.

For an appropriate subset $\{x_{t}:t\in T\}$ of $A_{0}$, where the indices are
taken from an idex set $T$, we have $\{\left[  x\right]  _{f^{2}}:x\in
A_{0}\}=\{\left[  x_{t}\right]  _{f^{2}}:t\in T\}$, and $\left[  x_{t}\right]
_{f^{2}}\neq\left[  x_{s}\right]  _{f^{2}}$ for all $t,s\in T$ with $t\neq s$.
Such a subset $\{x_{t}:t\in T\}\subseteq A_{0}$ is an \textit{irredundant set
of representatives} of the equivalence classes of $\sim_{f^{2}}$ (in $A_{0}$).
That is%
\[
A_{0}=%
{\textstyle\bigcup\nolimits_{t\in T}}
\left[  x_{t}\right]  _{f^{2}}\text{ and }\left[  x_{t}\right]  _{f^{2}}%
\cap\left[  x_{s}\right]  _{f^{2}}=\varnothing\text{ for all }t,s\in T\text{
with }t\neq s.
\]
Now for each $t\in T$ consider the extension $R_{t}=R\cap(\left[
x_{t}\right]  _{f^{2}}\times\left[  x_{t}\right]  _{f^{2}})$ of the restricted
partial order $r\cap(\left[  x_{t}\right]  _{f^{2}}\times\left[  x_{t}\right]
_{f^{2}})$. By Corollary 3.2, the relation
\[
R_{t}\cup R_{t}^{-1}=(\left[  x_{t}\right]  _{f^{2}}\times\left[
x_{t}\right]  _{f^{2}})\smallsetminus\pi=(\left[  x_{t}\right]  _{f^{2}}%
\times\left[  x_{t}\right]  _{f^{2}})\smallsetminus\beta
\]
of $R_{t}$-comparability is an equivalence. Thus $\left[  x_{t}\right]
_{f^{2}}$ is a disjoint union of finitely many $R_{t}$-chains (these are the
equivalence classes with respect to $(\left[  x_{t}\right]  _{f^{2}}%
\times\left[  x_{t}\right]  _{f^{2}})\smallsetminus\beta$) and any two
elements from different chains are incomparable with respect to $R_{t}$.

We claim that
\[
S=R_{\ast}\cup\left(
{\textstyle\bigcup\nolimits_{t\in T}}
R_{t}\right)  \cup P\cup Q
\]
with%
\[
P=\{(a,x):a\in A_{\ast},x\in A_{0}\text{ and }a\leq_{R_{\ast}}p\}
\]
and%
\[
Q=\{(y,b):b\in A_{\ast},y\in A_{0}\text{ and }q\leq_{R_{\ast}}b\}
\]
is a lattice order extension of $r$\ and that $f$ is a lattice
anti-endomorphism of $(A,\leq_{S})$.

The proof consists of the following straightforward steps.

Notice that, $P\subseteq A_{\ast}\times A_{0}$, and $Q\subseteq A_{0}\times
A_{\ast}$. Also the direct products $A_{\ast}\times A_{\ast}$, $A_{\ast}\times
A_{0}$, $A_{0}\times A_{\ast}$, and $\left[  x_{t}\right]  _{f^{2}}%
\times\left[  x_{t}\right]  _{f^{2}}$ (for $t\in T$) are pairwise disjoint.

In order to see $r\subseteq S$, take $(u,v)\in r$.

\noindent(1) If $(u,v)\in A_{\ast}\times A_{\ast}$, then $r\cap(A_{\ast}\times
A_{\ast})\subseteq R\cap(A_{\ast}\times A_{\ast})=R_{\ast}\subseteq S$ implies
$(u,v)\in R$.

\noindent(2) If $(u,v)\in A_{\ast}\times A_{0}$, then $\left[  u\right]
_{f^{2}}\neq\left[  v\right]  _{f^{2}}$ and $2\neq n_{f}(v)\neq\infty$
contradicts $(u,v)\in r$.

\noindent(3) $(u,v)\in A_{0}\times A_{\ast}$ is also impossible.

\noindent(4) If $(u,v)\in A_{0}\times A_{0}$, then $(u,v)\in\left[
x_{t}\right]  _{f^{2}}\times\left[  x_{s}\right]  _{f^{2}}$ for some $t,s\in
T$. Clearly, $t\neq s$ would imply $\left[  u\right]  _{f^{2}}\neq\left[
v\right]  _{f^{2}}$, and then $2\neq n_{f}(v)\neq\infty$ contradicts $(u,v)\in
r$. Thus $t=s$, and $r\cap(\left[  x_{t}\right]  _{f}\times\left[
x_{t}\right]  _{f})\subseteq R\cap(\left[  x_{t}\right]  _{f^{2}}\times\left[
x_{t}\right]  _{f^{2}})=R_{t}$ yields $(u,v)\in S$.

We prove that $S$ is a partial order.

\noindent Antisymmetry: Let $(u,v)\in S$ and $(v,u)\in S$.

\noindent(1) If $(u,v),(v,u)\in R_{\ast}$, then $u=v$ follows from the
antisymmetric property of $R_{\ast}$

\noindent(2) If $(u,v)\in R_{t}$ and $(v,u)\in R_{s}$, then $t=s$, and $u=v$
follows from the antisymmetric property of $R_{t}$.

\noindent(3) If $(u,v)\in P$ and $(v,u)\in Q$, then $u\leq_{R_{\ast}}p$ and
$q\leq_{R_{\ast}}u$ imply $q\leq_{R_{\ast}}p$, contradicting with
$p\leq_{R_{\ast}}q$ and $p\neq q$.

\noindent(4) If $(u,v)\in Q$ and $(v,u)\in P$, then interchanging the roles of
$u$ and $v$ leads to a similar contradiction as in case (3).

\noindent Transitivity: Let $(u,v)\in S$ and $(v,w)\in S$.

\noindent(1) If $(u,v),(v,w)\in R_{\ast}$, then $(u,w)\in R_{\ast}$ follows
from the transitivity of $R_{\ast}$.

\noindent(2) If $(u,v)\in R_{\ast}$ and $(v,w)\in P$, then $u\leq_{R_{\ast}%
}v\leq_{R_{\ast}}p$ and $w\in A_{0}$ imply $(u,w)\in P$.

\noindent(3) If $(u,v)\in R_{t}$ and $(v,w)\in R_{s}$, then we have $t=s$, and
$(u,w)\in R_{t}$ follows from the transitivity of $R_{t}$.

\noindent(4) If $(u,v)\in R_{t}$ and $(v,w)\in Q$, then $u,v\in A_{0}$, $w\in
A_{\ast}$, and $q\leq_{R_{\ast}}w$. It follows that $(u,w)\in Q$.

\noindent(5) If $(u,v)\in P$ and $(v,w)\in R_{t}$, then $v,w\in A_{0}$, $u\in
A_{\ast}$, and $u\leq_{R_{\ast}}p$. It follows that $(u,w)\in P$.

\noindent(6) If $(u,v)\in P$ and $(v,w)\in Q$, then $u\leq_{R_{\ast}}%
p\leq_{R_{\ast}}q\leq_{R_{\ast}}w$, from which $(u,w)\in R_{\ast}$ follows.

\noindent(7) If $(u,v)\in Q$ and $(v,w)\in R_{\ast}$, then $u\in A_{0}$ and
$q\leq_{R_{\ast}}v\leq_{R_{\ast}}w$ imply $(u,w)\in P$.

\noindent(8) If $(u,v)\in Q$ and $(v,w)\in P$, then $q\leq_{R_{\ast}}%
v\leq_{R_{\ast}}p$ contradicts $p\leq_{R_{\ast}}q$ and $p\neq q$.

We note that $f$\ is order-reversing with respect to $(A_{\ast},\leq_{R_{\ast
}})$, and $(\left[  x_{t}\right]  _{f^{2}},\leq_{R_{t}})$ for $t\in T$. In
order to check the order-reversing property of $f$ with respect to
$(A,\leq_{S})$, it is enough to see that $(a,x)\in P$ implies $(f(x),f(a))\in
Q$ and $(y,b)\in Q$ implies $(f(b),f(y))\in P$. Obviously, $a\in A_{\ast}$,
$x\in A_{0}$, and $a\leq_{R_{\ast}}p$ imply $f(a)\in A_{\ast}$, $f(x)\in
A_{0}$, and $q=f(p)\leq_{R_{\ast}}f(a)$. Similarly, $b\in A_{\ast}$, $y\in
A_{0}$, and $q\leq_{R_{\ast}}b$ imply $f(b)\in A_{\ast}$, $f(y)\in A_{0}$, and
$f(b)\leq_{R_{\ast}}f(q)=p$.

If $u,v\in A$ are comparable elements with respect to $S$, then the existence
of the supremum $u\vee v$ and the infimum $u\wedge v$ in $(A,\leq_{S})$ is
evident; moreover, the order-reversing property of $f$ ensures that%
\[
f(u\vee v)=f(u)\wedge f(v)\text{ and }f(u\wedge v)=f(u)\vee f(v).
\]

If $u,v\in A$ are incomparable elements with respect to $S$, then we have the
following possibilities.

\noindent(1) If $u\in A_{\ast}$ and $v\in A_{0}$, then $(u,v)\notin P$,
$(v,u)\notin Q$, and the linearity of $R_{\ast}$\ imply $p\leq_{R_{\ast}}%
u\leq_{R_{\ast}}q$, from which $u\vee v=q$ and $u\wedge v=p$ follow in
$(A,\leq_{S})$. Since $f(u)\in A_{\ast}$, $f(v)\in A_{0}$, and $p=f(q)\leq
_{R_{\ast}}f(u)\leq_{R_{\ast}}f(p)=q$, we deduce that%
\[
f(u\vee v)=f(q)=p=f(u)\wedge f(v)\text{ and }f(u\wedge v)=f(p)=q=f(u)\vee
f(v).
\]
\noindent(2) If $u\in A_{0}$ and $v\in A_{\ast}$, then interchanging the roles
of $u$ and $v$ leads to the same result as in case (1).

\noindent(3) If $u,v\in A_{0}$ and $\left[  u\right]  _{f^{2}}\neq\left[
v\right]  _{f^{2}}$, then $u\vee v=q$ and $u\wedge v=p$ in $(A,\leq_{S})$
follow directly from the definition of $S$. Since $f(u),f(v)\in A_{0}$ and
$\left[  u\right]  _{f^{2}}\neq\left[  v\right]  _{f^{2}}$ implies $\left[
f(u)\right]  _{f^{2}}\neq\left[  f(v)\right]  _{f^{2}}$, we deduce%
\[
f(u\vee v)=f(q)=p=f(u)\wedge f(v)\text{ and }f(u\wedge v)=f(p)=q=f(u)\vee
f(v).
\]
\noindent(4) If $u,v\in A_{0}$ and $\left[  u\right]  _{f^{2}}=\left[
v\right]  _{f^{2}}=\left[  x_{t}\right]  _{f^{2}}$ for some unique $t\in T$,
then $(u,v)\notin R_{t}$ and $(v,u)\notin R_{t}$ imply $(u,v)\in\beta$. It
follows that $u$ and $v$ are in different equivalence classes with respect to
the $R_{t}$-comparability relation $(\left[  x_{t}\right]  _{f^{2}}%
\times\left[  x_{t}\right]  _{f^{2}})\smallsetminus\beta$. An upper (lower)
bound of $\{u,v\}$ in $(\left[  x_{t}\right]  _{f^{2}},\leq_{R_{t}})$ would be
comparable with $u$ and $v$, which is impossible. We conclude that the set
$\{u,v\}$ has no upper and lower bounds in $(\left[  x_{t}\right]  _{f^{2}%
},\leq_{R_{t}})$. Thus we have $u\vee v=q$ and $u\wedge v=p$ in $(A,\leq_{S})$.

Since $f(u),f(v)\in\left[  f(x_{t})\right]  _{f^{2}}$ and $(f(u),f(v))\in
\beta$ (by Lemma 3.1), a similar argument (in $(\left[  x_{s}\right]  _{f^{2}%
},\leq_{R_{s}})$\ with $\left[  x_{s}\right]  _{f^{2}}=\left[  f(x_{t}%
)\right]  _{f^{2}}$) gives that%
\[
f(u)\vee f(v)=q=f(p)=f(u\wedge v)\text{ and }f(u)\wedge f(v)=p=f(q)=f(u\vee
v).\square
\]

\bigskip

\noindent Theorems 3.7 and 3.8 together generalize the answer given in section
2 to the question in the title of the paper. We pose a further problem.

\bigskip

\noindent\textbf{Problem.}\textit{ Consider an arbitrary function
}$f:A\longrightarrow A$\textit{. Find necessary and sufficient conditions for
the existence of a modular (or distributive) lattice structure }%
$(A,\vee,\wedge)$\textit{ on }$A$\textit{ such that }$f$\textit{ becomes a
lattice anti-endomorphism of }$(A,\vee,\wedge)$\textit{.}

\bigskip

\noindent\textbf{Example.} Let $A=\{p,q,x_{1},x_{2},\ldots,x_{n}\}$, where
$n=2k+1\geq3$ is odd, and let $f:A\longrightarrow A$ be a function with
$f(p)=q$, $f(q)=p$, $f(x_{n})=x_{1}$, and $f(x_{i})=x_{i+1}$ for $1\leq i\leq
n-1$. If $f$ is an anti-endomorphism of some lattice $(A,\leq,\vee,\wedge)$,
then $f$ is order-reversing with respect to $(A,\leq)$, and Proposition 3.3
ensures that the proper cycle $\{x_{1},\ldots,x_{n}\}$ of $f^{2}$\ is an
antichain in $(A,\leq)$. Since $x_{1}\vee\cdots\vee x_{n}$ and $x_{1}%
\wedge\cdots\wedge x_{n}$ form a two-element cycle of $f$, one of $x_{1}%
\vee\cdots\vee x_{n}$ and $x_{1}\wedge\cdots\wedge x_{n}$ is $p$ and the other
is $q$. Thus $(A,\leq,\vee,\wedge)$ is isomorphic to the lattice $M_{n}$ in
both cases. It follows that there is no distributive lattice structure on $A$
making $f$ a lattice anti-endomorphism (even though $f$ has a cycle of length
$2$).

\bigskip

\noindent REFERENCES

\bigskip

\noindent\lbrack FSz] Foldes, S.; Szigeti, J. \textit{Maximal compatible
extensions of partial orders, }J. Australian Math. Soc. 81 (2006), 245-252.

\noindent\lbrack G] Gr\"{a}tzer, G. \textit{General Lattice Theory,}
Birkhauser Verlag, Basel-Boston-Berlin (2003)

\noindent\lbrack JPR] Jakub\'{\i}kova-Studenovsk\'{a}, D., P\"{o}schel, R. and
Radeleczki, S.: \textit{The lattice of compatible quasiorders of acyclic
monounary algebras,} Order 28 (2011), 481-497.

\noindent\lbrack L] Lengv\'{a}rszky, Zs. \textit{Linear extensions of partial
orders preserving antimonotonicity,} Publicationes Math. Debrecen 38 (1991),
no.3-4, 279-285.

\noindent\lbrack Sz1] Szigeti, J. \textit{Maximal extensions of partial orders
preserving antimonotonicity}, Algebra Universalis, Vol. 66, No. 1-2 (2011), 143-150.

\noindent\lbrack Sz2] Szigeti, J. \textit{Which self-maps appear as lattice
endomorphisms?, }Discrete Mathematics, 321 (2014), 53--56.

\end{document}